\documentclass[12pt]{article}
\usepackage{wasysym,amssymb,eufrak,indentfirst,cite,bibspacing,amsthm}
\begin{document}
\baselineskip 17pt
\title{\textbf{On $\Pi$-supplemented subgroups of a finite group}\thanks{Research is supported by a NNSF grant of China (grant \#11371335) and Research Fund for the Doctoral Program of Higher Education of China (Grant 20113402110036).}}
\author{Xiaoyu Chen, Wenbin Guo\thanks{Corresponding author.}\\
{\small Department of Mathematics, University of Science and Technology of China,}\\ {\small Hefei 230026, P. R. China}\\
 {\small E-mail: jelly@mail.ustc.edu.cn, $\,$wbguo@ustc.edu.cn}}
\date{}
\maketitle
\begin{abstract}
A subgroup $H$ of a finite group $G$ is said to satisfy $\Pi$-property in $G$ if for every chief factor $L/K$ of $G$, $|G/K:N_{G/K}(HK/K\cap L/K)|$ is a $\pi(HK/K\cap L/K)$-number. A subgroup $H$ of $G$ is called to be $\Pi$-supplemented in $G$ if there exists a subgroup $T$ of $G$ such that $G=HT$ and $H\cap T\leq I\leq H$, where $I$ satisfies $\Pi$-property in $G$. In this paper, we investigate the structure of a finite group $G$ under the assumption that some primary subgroups of $G$ are $\Pi$-supplemented in $G$. The main result we proved improves a large number of earlier results.
\end{abstract}
\renewcommand{\thefootnote}{\empty}

\section{Introduction}
\noindent Throughout this paper, all groups mentioned are finite, $G$ always denotes a finite group and $p$ denotes a prime. Let $\pi$ denote a set of some primes, $\pi(G)$ denote the set of all prime divisors of $|G|$, and $|G|_p$ denote the order of the Sylow $p$-subgroups of $G$. An integer $n$ is called a $\pi$-number if all prime divisors of $n$ belong to $\pi$. For a subgroup $H$ of $G$, let $H^G$ denote the normal closure of $H$ in $G$, that is, $H^G=\langle H^g\mid g\in G\rangle$.\par
Recall that a class of groups $\mathcal{F}$ is called a formation if $\mathcal{F}$ is closed under taking homomorphic images and subdirect products. A formation $\mathcal{F}$ is said to be saturated (resp. solubly saturated) if $G\in \mathcal{F}$ whenever $G/\Phi(G)\in \mathcal{F}$ (resp. $G/\Phi(N)\in \mathcal{F}$ for a soluble normal subgroup $N$ of $G$). A chief factor $L/K$ of $G$ is said to be $\mathcal{F}$-central in $G$ if $(L/K)\rtimes (G/C_G(L/K))\in \mathcal{F}$. A normal subgroup $N$ of $G$ is called to be $\mathcal{F}$-hypercentral in $G$ if every chief factor of $G$ below $N$ is $\mathcal{F}$-central in $G$. Let $Z_\mathcal{F}(G)$ denote the $\mathcal{F}$-hypercentre of $G$, that is, the product of all $\mathcal{F}$-hypercentral normal subgroups of $G$. We use $\mathcal{U}$ (resp. $\mathcal{U}_p$) to denote the class of finite supersoluble (resp. $p$-supersoluble) groups and $\mathcal{G}_\pi$ to denote the class of all finite $\pi$-groups.\par
Recall that $G$ is said to be quasinilpotent if for every chief factor $L/K$ of $G$ and every element $x\in G$, $x$ induces an inner automorphism on $L/K$. The generalized Fitting subgroup $F^*(G)$ of $G$ is the quasinilpotent radical of $G$ (for details, see \cite[Chapter X]{Hup}). All notations and terminology not mentioned above are standard, as in \cite{Hup1,Doe,Guo1}.\par
In \cite{Li1}, Li introduced the concepts of $\Pi$-property and $\Pi$-supplemented subgroup as follows:\par
\medskip
\noindent\textbf{Definition 1.1.} \cite{Li1} A subgroup $H$ of $G$ is said to satisfy \textit{$\Pi$-property} in $G$ if for every chief factor $L/K$ of $G$, $|G/K:N_{G/K}(HK/K\cap L/K)|$ is a $\pi(HK/K\cap L/K)$-number.\par
A subgroup $H$ of $G$ is called to be \textit{$\Pi$-supplemented} in $G$ if there exists a subgroup $T$ of $G$ such that $G=HT$ and $H\cap T\leq I\leq H$, where $I$ satisfies $\Pi$-property in $G$.\par
\medskip
As we showed in Section 4 below, the concept of $\Pi$-supplemented subgroup generalizes many known embedding properties. However, besides \cite{Li6}, this concept has not been deeply investigated. In this paper, we will continue to study the properties of $\Pi$-supplemented subgroups, and arrive at the following main result.\par
\medskip
\noindent\textbf{Theorem A.} \textit{Let $\mathcal{F}$ be a solubly saturated formation containing $\mathcal{U}$ and $E$ a normal subgroup of $G$ with $G/E\in \mathcal{F}$. Let $X\unlhd G$ such that $F^*(E)\leq X\leq E$. For every prime $p\in\pi(X)$ and every non-cyclic Sylow $p$-subgroup $P$ of $X$, suppose that $P$ has a subgroup $D$ such that $1\leq |D|<|P|$ and every proper subgroup $H$ of $P$ with $|H|=p^n|D|$ $($$n=0,1$$)$ either is $\Pi$-supplemented in $G$ or has a $p$-supersoluble supplement in $G$. If $P$ is not quaternion-free and $|D|=1$, suppose further that every cyclic subgroup of $P$ of order $4$ either is $\Pi$-supplemented in $G$ or has a $2$-supersoluble supplement in $G$. Then $G\in \mathcal{F}$.}\par
\medskip
Recall that a subgroup $H$ of $G$ is said to be \textit{$c$-supplemented} \cite{Bal2} in $G$ if there exists a subgroup $T$ of $G$ such that $G=HT$ and $H\cap T\leq H_G$, where $H_G$ denotes the largest normal subgroup of $G$ contained in $H$. It is easy to find that all $c$-supplemented subgroups of $G$ are $\Pi$-supplemented in $G$, and the converse does not hold. For example, let $G=\langle a,b\,|\,a^5=b^4=1, b^{-1}ab=a^2\rangle$ and $H=\langle b^2 \rangle$. Then $H$ is $\Pi$-supplemented, but not $c$-supplemented in $G$. In \cite{Asa2}, M. Asaad proved the following excellent theorem.\par
\medskip
\noindent\textbf{Theorem 1.2.} \cite[Theorems 1.5 and 1.6]{Asa2} \textit{Let $\mathcal{F}$ be a saturated formation containing $\mathcal{U}$ and $E$ a normal subgroup of $G$ with $G/E\in \mathcal{F}$. Let $X\unlhd G$ such that $X=E$ or $X=F^*(E)$. For any Sylow subgroup $P$ of $X$, let $D$ be a subgroup of $P$ such that $1\leq |D|<|P|$. Suppose that every subgroup $H$ of $P$ with $|H|=p^n|D|$ $($$n=0,1$$)$ is $c$-supplemented in $G$. If $P$ is a non-abelian $2$-group and $|D|=1$, suppose further that every cyclic subgroup of $P$ of order $4$ is $c$-supplemented in $G$. Then $G\in \mathcal{F}$.}\par
\medskip
One can see that Theorem A can be viewed as a large improvement of M. Asaad's result. The following theorems are the main stages of the proof of Theorem A.\par
\medskip
\noindent\textbf{Theorem B.} \textit{Let $P$ be a normal $p$-subgroup of $G$. Suppose that $P$ has a subgroup $D$ such that $1\leq |D|<|P|$ and every proper subgroup $H$ of $P$ with $|H|=p^n|D|$ $($$n=0,1$$)$ either is $\Pi$-supplemented in $G$ or has a $p$-supersoluble supplement in $G$. If $P$ is not quaternion-free and $|D|=1$, suppose further that every cyclic subgroup of $P$ of order $4$ either is $\Pi$-supplemented in $G$ or has a $2$-supersoluble supplement in $G$. Then $P\leq Z_\mathcal{U}(G)$.}\par
\medskip
\noindent\textbf{Theorem C.} \textit{Let $E$ be a normal subgroup of $G$ and $P$ a Sylow $p$-subgroup of $E$ with $(|E|,p-1)=1$. Suppose that $P$ has a subgroup $D$ such that $1\leq |D|<|P|$ and every proper subgroup $H$ of $P$ with $|H|=p^n|D|$ $($$n=0,1$$)$ either is $\Pi$-supplemented in $G$ or has a $p$-supersoluble supplement in $G$. If $P$ is not quaternion-free and $|D|=1$, suppose further that every cyclic subgroup of $P$ of order $4$ either is $\Pi$-supplemented in $G$ or has a $2$-supersoluble supplement in $G$. Then $E$ is $p$-nilpotent.}\par
\medskip
Finally, the following corollaries can be deduced immediately from Theorem A.\par
\medskip
\noindent\textbf{Corollary 1.3.} \textit{Let $\mathcal{F}$ be a solubly saturated formation containing $\mathcal{U}$ and $E$ a normal subgroup of $G$ with $G/E\in \mathcal{F}$. Let $X\unlhd G$ such that $F^*(E)\leq X\leq E$. For every prime $p\in\pi(X)$ and every non-cyclic Sylow $p$-subgroup $P$ of $X$, suppose that every maximal subgroup of $P$ either is $\Pi$-supplemented in $G$ or has a $p$-supersoluble supplement in $G$. Then $G\in \mathcal{F}$.}\par
\medskip
\noindent\textbf{Corollary 1.4.} \textit{Let $\mathcal{F}$ be a solubly saturated formation containing $\mathcal{U}$ and $E$ a normal subgroup of $G$ with $G/E\in \mathcal{F}$. Let $X\unlhd G$ such that $F^*(E)\leq X\leq E$. For every prime $p\in\pi(X)$ and every non-cyclic Sylow $p$-subgroup $P$ of $X$, suppose that every cyclic subgroup of $P$ of prime order or order $4$ $($when $P$ is not quaternion-free$)$ either is $\Pi$-supplemented in $G$ or has a $p$-supersoluble supplement in $G$. Then $G\in \mathcal{F}$.}\par
\section{Basic Properties}
\noindent\textbf{Lemma 2.1.} \cite[Proposition 2.1]{Li1} \textit{Let $H\leq G$ and $N\unlhd G$.}\par
(1) \textit{If $H$ satisfies $\Pi$-property in $G$, then $HN/N$ satisfies $\Pi$-property in $G/N$.}\par
(2) \textit{If either $N\leq H$ or $(|H|,|N|)=1$ and $H$ is $\Pi$-supplemented in $G$, then $HN/N$ is $\Pi$-supplemented in $G/N$.}\par
\medskip
\noindent\textbf{Lemma 2.2.} \textit{Let $H\leq G$, $N\unlhd G$ such that $N\leq H$ and $P$ be a Sylow $p$-subgroup of $H$. Suppose that $P$ has a subgroup $D$ such that $|N|_p\leq |D|<|P|$ and every subgroup of $P$ of order $|D|$ either is $\Pi$-supplemented in $G$ or has a $p$-supersoluble supplement in $G$. Then every subgroup of $PN/N$ of order $|D|/|N|_p$ either is $\Pi$-supplemented in $G/N$ or has a $p$-supersoluble supplement in $G/N$.}\par
\smallskip
\noindent\textit{Proof.} Let $X/N$ be a subgroup of $PN/N$ of order $|D|/|N|_p$. Then $X=(P\cap X)N$, and so $X/N\cong P\cap X/P\cap N$. Hence $|P\cap X|=|D|$. By hypothesis, $P\cap X$ either is $\Pi$-supplemented in $G$ or has a $p$-supersoluble supplement in $G$. Then $P\cap X$ has a supplement $T$ in $G$ such that either $T$ is $p$-supersoluble or $P\cap X\cap T\leq I\leq P\cap X$, where $I$ satisfies $\Pi$-property in $G$. Obviously, $G/N=(X/N)(TN/N)$. Since $(|N:P\cap N|,|N:T\cap N|)=1$, $N=(P\cap N)(T\cap N)$. This deduces that $X\cap TN=(P\cap X)N\cap TN=(P\cap X\cap T)N$. Therefore, either $TN/N$ is $p$-supersoluble or $X/N\cap TN/N=(P\cap X\cap T)N/N\leq IN/N\leq X/N$, where $IN/N$ satisfies $\Pi$-property in $G/N$ by Lemma 2.1(1). Consequently, $X/N$ either is $\Pi$-supplemented in $G/N$ or has a $p$-supersoluble supplement in $G/N$.\qed\par
\medskip
For any function $f$: $\mathbb{P}\:\cup \:\{0\}\longrightarrow \mbox{\{formations of groups\}}$. Following \cite{Ski1}, let $$CF(f)=\{G \mbox{ is a group } | \mbox{ } G/C_G(H/K)\in f(0) \mbox{ for each non-abelian}\mbox{ chief factor }$$ $$H/K\mbox{ of } G
\mbox{ and } G/C_G(H/K)\in f(p) \mbox{ for each abelian } p\mbox{-chief factor } H/K \mbox{ of } G\}\mbox{.}$$\par
\medskip
\noindent\textbf{Lemma 2.3.} \cite{Ski1} \textit{For any non-empty solubly saturated formation $\mathcal{F}$, there exists a unique function $F$: $\mathbb{P}\:\cup \:\{0\}\longrightarrow \mbox{\{formations of groups\}}$ such that $\mathcal{F} = CF(F)$, $F(p)=\mathcal{G}_pF(p)\subseteq \mathcal{F}$ for all $p\in\mathbb{P}$ and $F(0)=\mathcal{F}$.}\par
\medskip
The function $F$ in Lemma 2.3 is called the canonical composition satellite of $\mathcal{F}$.\par
\medskip
\noindent\textbf{Lemma 2.4.} \cite[Lemma 2.14]{Guo} \textit{Let $\mathcal{F}$ be a saturated $($resp. solubly saturated$)$ formation and $F$ the canonical local $($resp. the canonical composition$)$ satellite of $\mathcal{F}$ $($for the details of canonical local satellite, see \textup{\cite[\textit{Chapter IV, Definition $3.9$}]{Doe}}$)$. Let $E$ be a normal $p$-subgroup of $G$. Then $E\leq Z_\mathcal{F}(G)$ if and only if $G/C_G(E)\in F(p)$.}\par
\medskip
\noindent\textbf{Lemma 2.5.} \cite[Lemma 2.4]{Gag} \textit{Let $P$ be a $p$-group. If $\alpha$ is a $p'$-automorphism of $P$ which centralizes ${\Omega}_1(P)$, then $\alpha=1$ unless $P$ is a non-abelian $2$-group. If $[\alpha, {\Omega}_2(P)] = 1$, then $\alpha=1$ without restriction.}\par
\medskip
\noindent\textbf{Lemma 2.6.} \cite[Lemma 2.15]{Dor} \textit{If $\sigma$ is an automorphism of odd order of the quaternion-free $2$-group $P$ and $\sigma$ acts trivially  on ${\Omega}_1(P)$, then $\sigma=1$.}\par
\medskip
\noindent\textbf{Lemma 2.7.} \cite[Chapter 5, Theorem 3.13]{Gor} \textit{A $p$-group $P$ possesses a characteristic subgroup $C$ $($which is called a Thompson critical subgroup of $P$$)$ with the following properties:}\par
(1) \textit{The nilpotent class of $C$ is at most $2$, and $C/Z(C)$ is elementary abelian.}\par
(2) \textit{$[P,C]\leq Z(C)$.}\par
(3) \textit{$C_P(C)=Z(C)$.}\par
(4) \textit{Every nontrivial $p'$-automorphism of $P$ induces a nontrivial automorphism of $C$.}\par
\medskip
If $P$ is either an odd order $p$-group or a quaternion-free 2-group, then let ${\Omega}(P)$ denote the subgroup ${\Omega}_1(P)$, otherwise ${\Omega}(P)$ denotes ${\Omega}_2(P)$. The following lemma is a generalization of \cite[Lemma 2.12]{Che2}, which is attributed to A. N. Skiba.\par
\medskip
\noindent\textbf{Lemma 2.8.} \textit{Let $\mathcal{F}$ be a solubly saturated formation, $P$ a normal $p$-subgroup of $G$ and $C$ a Thompson critical subgroup of $P$. If either $P/\Phi(P)\leq Z_\mathcal{F}(G/\Phi(P))$ or ${\Omega}(C)\leq Z_\mathcal{F}(G)$, then $P\leq Z_\mathcal{F}(G)$.}\par
\smallskip
\noindent\textit{Proof.} Let $F$ be the canonical composition satellite of $\mathcal{F}$. Suppose that $P/\Phi(P)\leq Z_\mathcal{F}(G/\Phi(P))$. Then by Lemma 2.4, $G/C_G(P/\Phi(P))\in F(p)$. Note that by \cite[Chapter 5, Theorem 1.4]{Gor}, $C_G(P/\Phi(P))/C_G(P)$ is a $p$-group. This implies that $G/C_G(P)\in \mathcal{G}_pF(p)=F(p)$. Hence by Lemma 2.4 again, $P\leq Z_\mathcal{F}(G)$.\par
Now assume that ${\Omega}(C)\leq Z_\mathcal{F}(G)$. Then by Lemma 2.4, $G/C_G({\Omega}(C))\in F(p)$. Since $C_G({\Omega}(C))/C_G(C)$ is a $p$-group by Lemmas 2.5 and 2.6, we have that $G/C_G(C)\in \mathcal{G}_pF(p)=F(p)$. It follows from Lemma 2.7(4) that $C_G(C)/C_G(P)$ is a $p$-group, and so $G/C_G(P)\in \mathcal{G}_pF(p)=F(p)$. Thus by Lemma 2.4 again, $P\leq Z_\mathcal{F}(G)$.\qed\par
\medskip
\noindent\textbf{Lemma 2.9.} \cite[Lemma 3.1]{War} \textit{Let $G$ be a non-abelian quaternion-free $2$-group. Then $G$ has a characteristic subgroup of index $2$.}\par
\medskip
\noindent\textbf{Lemma 2.10.} \textit{Let $C$ be a Thompson critical subgroup of a nontrivial $p$-group $P$.}\par
(1) \textit{If $p$ is odd, then the exponent of ${\Omega}_1(C)$ is $p$.}\par
(2) \textit{If $P$ is an abelian $2$-group, then the exponent of ${\Omega}_1(C)$ is $2$.}\par
(3) \textit{If $p=2$, then the exponent of ${\Omega}_2(C)$ is at most $4$.}\par
\smallskip
\noindent\textit{Proof.} (1) Since the nilpotent class of $C$ is at most 2 by Lemma 2.7(1), the statement (1) directly follows from \cite[Chapter 5, Lemma 3.9(i)]{Gor}.\par
Statement (2) is clear.\par
(3) Let $x$ and $y$ be elements of $C$ of order 4. Then by Lemma 2.7(1) and \cite[Chapter 2, Lemma 2.2]{Gor}, we have that $[x,y]^2=[x^2,y]=1$, and so $(yx)^4=[x,y]^6y^4x^4=1$. This shows that the order of $yx$ is at most 4, and thus the exponent of ${\Omega}_2(C)$ is at most 4.\qed\par
\medskip
Recall that $G$ is said to be $\pi$-closed if $G$ has a normal Hall $\pi$-subgroup. Also, $G$ is said to be a $C_\pi$-group if $G$ has a Hall $\pi$-subgroup and any two Hall $\pi$-subgroups of $G$ are conjugate in $G$.\par
\medskip
\noindent\textbf{Lemma 2.11.} \textit{Let $p$ be a prime divisor of $|G|$ with $(|G|,p-1)=1$. Suppose that $G$ has a Hall $p'$-subgroup. Then $G$ is a $C_{p'}$-group.}\par
\smallskip
\noindent\textit{Proof.} If $p>2$, then $2\nmid |G|$. By Feit-Thompson Theorem, $G$ is soluble, and so $G$ is a $C_{p'}$-group. If $p=2$, then by \cite[Theorem A]{Cro}, $G$ is also a $C_{p'}$-group.\qed\par
\medskip
\noindent\textbf{Lemma 2.12.} \cite[Corollary 3.7]{Guo} \textit{Let $P$ be a $p$-subgroup of $G$. Suppose that $G$ is a $C_\pi$-group for some set of primes $\pi$ with $p\notin \pi$. If every maximal subgroup of $P$ has a $\pi$-closed supplement in $G$, then $G$ is $\pi$-closed.}\par
\medskip
The next lemma is well-known.\par
\medskip
\noindent\textbf{Lemma 2.13.} \textit{Let $p$ be a prime divisor of $|G|$ with $(|G|,p-1)=1$.}\par
(1) \textit{If $G$ has cyclic Sylow $p$-subgroups, then $G$ is $p$-nilpotent.}\par
(2) \textit{If $E$ is a normal subgroup of $G$ such that $|E|_p\leq p$ and $G/E$ is $p$-nilpotent, then $G$ is $p$-nilpotent.}\par
(3) \textit{If $H$ is a subgroup of $G$ such that $|G:H|=p$, then $H\unlhd G$.}\par
\medskip
\noindent\textbf{Lemma 2.14.} \cite[Lemma 2.6]{Asa1} \textit{If $G$ possesses two subgroups $K$ and $T$ such that $|G:K|=2^r$ and $|G:T|=2^{r+1}$ $(r\geq 3)$ and $T$ is not a $2'$-Hall subgroup of $G$, then $G$ is not a non-abelian simple group.}\par
\medskip
Recall that a subgroup $H$ of $G$ is said to be \textit{complemented} in $G$ if there exists a subgroup $T$ of $G$ such that $G=HT$ and $H\cap T=1$. In this case, $T$ is called a complement of $H$ in $G$.\par
\medskip
\noindent\textbf{Lemma 2.15.} \textit{Let $P$ be a Sylow $p$-subgroup of $G$ with $(|G|,p-1)=1$. If every subgroup of $P$ of order $p$ is complemented in $G$, then $G$ is $p$-nilpotent.}\par
\smallskip
\noindent\textit{Proof.} Let $H$ be a subgroup of $P$ of order $p$ and $T$ a complement of $H$ in $G$. Then by Lemma 2.13(3), $T\unlhd G$. If $p\nmid |T|$, then $G$ is $p$-nilpotent. Thus $p\mid |T|$. Clearly, $P\cap T$ is a Sylow $p$-subgroup of $T$ and every subgroup of $P\cap T$ of order $p$ is complemented in $T$. Then by induction, $T$ is $p$-nilpotent. Since the normal $p$-complement of $T$ is the normal $p$-complement of $G$, $G$ is also $p$-nilpotent.\qed\par
\medskip
\noindent\textbf{Lemma 2.16.} \cite[Lemma 3.1]{Bal3} \textit{Let $\mathfrak{F}$ be a saturated formation of characteristic $\pi$ and $H$ a subnormal subgroup of $G$ containing $O_\pi(\Phi(G))$ such that $H/O_\pi(\Phi(G))\in \mathfrak{F}$. Then $H\in \mathfrak{F}$.}\par
\medskip
\noindent\textbf{Lemma 2.17.} \cite[Theorem B]{Ski} \textit{Let $\mathcal{F}$ be any formation. If $E\unlhd G$ and $F^*(E)\leq Z_\mathcal{F}(G)$, then $E\leq Z_\mathcal{F}(G)$.}\par
\medskip
\noindent\textbf{Lemma 2.18.} \cite[Lemma 2.13]{Guo} \textit{Let $\mathcal{F}=CF(F)$ be a solubly saturated formation, where $F$ is the canonical composition satellite of $\mathcal{F}$. Let $H/K$ be a chief factor of $G$. Then $H/K$ is $\mathcal{F}$-central in $G$ if and only if $G/C_G(H/K)\in F(p)$ in the case where $H/K$ is a $p$-group, and $G/C_G(H/K)\in F(0)=\mathcal{F}$ in the case where $H/K$ is non-abelian.}\par
\section{Proofs of Theorems}
\noindent\textbf{Proof of Theorem B.} Suppose that the result is false and let $(G,P)$ be a counterexample for which $|G|+|P|$ is minimal. We proceed via the following steps.\par
(1) \textit{$|D|\geq p^2$.}\par
If $|D|\leq p$, we may assume that $|D|=1$ (in the conditions of the theorem, the case $|D|=p$ can be viewed as a special case of $|D|=1$). Then:\par
(i) \textit{$G$ has a unique normal subgroup $N$ such that $P/N$ is a chief factor of $G$, $N\leq Z_\mathcal{U}(G)$ and $|P/N|>p$.}\par
Let $P/N$ be a chief factor of $G$. Then $(G,N)$ satisfies the hypothesis of this theorem. By the choice of $(G,P)$, we have that $N\leq Z_\mathcal{U}(G)$. If $P/N\leq Z_\mathcal{U}(G/N)$, then $P\leq Z_\mathcal{U}(G)$, which is impossible. Hence $P/N\nleq Z_\mathcal{U}(G/N)$, and so $|P/N|>p$. Now let $P/R$ be a chief factor of $G$, which is different from $P/N$. Then we can obtain that $R\leq Z_\mathcal{U}(G)$ similarly as above. This implies that $P/N\leq Z_\mathcal{U}(G/N)$ by $G$-isomorphism $P/N=NR/N\cong R/N\cap R$, a contradiction.\par
(ii) \textit{Let $C$ be a Thompson critical subgroup of $P$. Then $P={\Omega}(C)$.}\par
If not, then ${\Omega}(C)\leq N\leq Z_\mathcal{U}(G)$ by (i). Thus by Lemma 2.8, $P\leq Z_\mathcal{U}(G)$, which is absurd.\par
(iii) \textit{The exponent of $P$ is $p$ or $4$ $($when $P$ is not quaternion-free$)$.}\par
If $P$ is a non-abelian quaternion-free 2-group, then $P$ has a characteristic subgroup $T$ of index 2 by Lemma 2.9. It follows from (i) that $T\leq N$, and so $|P/N|=2$, which is impossible. Hence by (ii) and Lemma 2.10, the exponent of $P$ is $p$ or 4 (when $P$ is not quaternion-free).\par
(iv) \textit{Final contradiction of $(1)$.}\par
Let $G_p$ be a Sylow $p$-subgroup of $G$. Since $P/N\cap Z(G_p/N)>1$, we may take a subgroup $V/N$ of $P/N\cap Z(G_p/N)$ of order $p$. Let $l\in V\backslash N$ and $H=\langle l\rangle$. Then $V=HN$ and $H$ is a group of order $p$ or 4 (when $P$ is not quaternion-free) by (iii). By hypothesis, $H$ either is $\Pi$-supplemented in $G$ or has a $p$-supersoluble supplement in $G$. Let $X$ be any supplement of $H$ in $G$. If $P\nleq X$, then $P\cap X<P$. Since $(P\cap X)^G=(P\cap X)^P<P$, we have that $P\cap X\leq N$ by (i). This implies that $P/N$ is cyclic for $P/P\cap X\cong H/H\cap X$ is cyclic, and so $|P/N|=p$, which  contradicts (i). Therefore, $P\leq X$, and thereby $X=G$. Consequently, $G$ is the unique supplement of $H$ in $G$. If $H$ has a $p$-supersoluble supplement in $G$, then $G$ is $p$-supersoluble. It follows that $P\leq Z_\mathcal{U}(G)$, which is impossible. Hence $H$ is $\Pi$-supplemented in $G$, and so $H$ satisfies $\Pi$-property in $G$. Then $|G:N_G(V)|=|G:N_G(HN)|$ is a $p$-number. This induces that $V\unlhd G$. Then by (i), $P=V$, and so $|P/N|=p$, a contradiction. This completes the proof of (1).\par
(2) \textit{$\Phi(P)=1$, and so $P$ is an elementary abelian $p$-group.}\par
Suppose that $\Phi(P)>1$. If $|\Phi(P)|>|D|$, then $(G,\Phi(P))$ satisfies the hypothesis of this theorem. By the choice of $(G,P)$, we have that $\Phi(P)\leq Z_\mathcal{U}(G)$. Let $L$ be a minimal normal subgroup of $G$ contained in $\Phi(P)$. Then $|L|=p$. Since $|D|>|L|=p$ by (1), $(G/L, P/L)$ satisfies the hypothesis of this theorem by Lemma 2.1(2). By the choice of $(G,P)$, $P/L\leq Z_\mathcal{U}(G/L)$. It follows that $P\leq Z_\mathcal{U}(G)$, which is absurd.\par
Hence $|\Phi(P)|\leq |D|$. Now we shall show that $P/\Phi(P)\leq Z_\mathcal{U}(G/\Phi(P))$. If $|\Phi(P)|<|D|$, then by Lemma 2.1(2), $(G/\Phi(P), P/\Phi(P))$ satisfies the hypothesis of this theorem. The choice of $(G,P)$ implies that $P/\Phi(P)\leq Z_\mathcal{U}(G/\Phi(P))$. Hence we may consider that $|\Phi(P)|=|D|$. If $p|D|=|P|$, then clearly, $P/\Phi(P)\leq Z_\mathcal{U}(G/\Phi(P))$. If $p|D|<|P|$, then by Lemma 2.1(2), every subgroup of $P/\Phi(P)$ of order $p$ either is $\Pi$-supplemented in $G/\Phi(P)$ or has a $p$-supersoluble supplement in $G/\Phi(P)$. This shows that $(G/\Phi(P), P/\Phi(P))$ satisfies the hypothesis of this theorem, and so $P/\Phi(P)\leq Z_\mathcal{U}(G/\Phi(P))$ by the choice of $(G,P)$. Then by Lemma 2.8, $P\leq Z_\mathcal{U}(G)$, which is impossible. Therefore, $\Phi(P)=1$.\par
(3) \textit{$G$ has a unique minimal normal subgroup $N$ contained in $P$, $P/N\leq Z_\mathcal{U}(G/N)$ and $p<|N|\leq |D|$.}\par
Let $G_p$ be a Sylow $p$-subgroup of $G$ and $N$ a minimal normal subgroup of $G$ contained in $P$. If $N=P$, then $P$ is a minimal normal subgroup of $G$. Let $H$ be a subgroup of $P$ of order $|D|$ such that $H\unlhd G_p$. By hypothesis, $H$ either is $\Pi$-supplemented in $G$ or has a $p$-supersoluble supplement in $G$. For any supplement $X$ of $H$ in $G$, we have that $P\cap X\unlhd G$. If $P\cap X=1$, then $H=P$, which is impossible. This induces that $P\cap X=P$, and so $X=G$. Therefore, $G$ is the unique supplement of $H$ in $G$. Since $G$ is not $p$-supersoluble, $H$ satisfies $\Pi$-property in $G$. It follows that $|G:N_G(H)|$ is a $p$-number. Hence $H\unlhd G$, a contradiction. Consequently, $N<P$. If $|N|>|D|$, then $(G,N)$ satisfies the hypothesis of this theorem. By the choice of $(G,P)$, we have that $N\leq Z_\mathcal{U}(G)$. This shows that $|N|=p>|D|$, which contradicts (1). Therefore, $|N|\leq |D|$.\par
Now we claim that $P/N\leq Z_\mathcal{U}(G/N)$. If $|N|<|D|$, then by Lemma 2.1(2), $(G/N, P/N)$ satisfies the hypothesis of this theorem. By the choice of $(G,P)$, $P/N\leq Z_\mathcal{U}(G/N)$. Hence we may assume that $|N|=|D|$. If $p|D|=|P|$, then clearly, $P/N\leq Z_\mathcal{U}(G/N)$. If $p|D|<|P|$, then by Lemma 2.1(2), every subgroup of $P/N$ of order $p$ either is $\Pi$-supplemented in $G/N$ or has a $p$-supersoluble supplement in $G/N$. Since $P$ is abelian, $(G/N,P/N)$ satisfies the hypothesis of this theorem. Then by the choice of $(G,P)$, we also have that $P/N\leq Z_\mathcal{U}(G/N)$. Consequently, our claim holds. If $|N|=p$, then $N\leq Z_\mathcal{U}(G)$, and so $P\leq Z_\mathcal{U}(G)$, which is absurd. Thus $|N|>p$. If $G$ has a minimal normal subgroup $R$ contained in $P$, which is different from $N$, then we get that $G/R\leq Z_\mathcal{U}(G/R)$ similarly as above. It follows that $NR/R\leq Z_\mathcal{U}(G/R)$, and so $N\leq Z_\mathcal{U}(G)$ for $G$-isomorphism $N\cong NR/R$. This implies that $P\leq Z_\mathcal{U}(G)$, a contradiction. Hence (3) holds.\par
(4) \textit{$p|D|=|P|$.}\par
If $p|D|<|P|$, then since $P/N\leq Z_\mathcal{U}(G/N)$, $G$ has a normal subgroup $K$ properly contained in $P$ such that $N\leq K$ and $|K|=p|D|$. Then $(G,K)$ satisfies the hypothesis of this theorem. By the choice of $(G,P)$, we have that $K\leq Z_\mathcal{U}(G)$, and thus $|N|=p$, which contradicts (3). This shows that (4) holds.\par
(5) \textit{Final contradiction.}\par
Since $\Phi(P)=1$, $N$ has a complement $S$ in $P$. Let $L$ be a maximal subgroup of $N$ such that $L\unlhd G_p$. Then $L\neq 1$ and $H=LS$ is a maximal subgroup of $P$. By hypothesis and (4), $H$ either is $\Pi$-supplemented in $G$ or has a $p$-supersoluble supplement in $G$. For any supplement $X$ of $H$ in $G$, since $P$ is abelian, we have that $P\cap X\unlhd G$. If $P\cap X=1$, then $H=P$, which is impossible. Hence $P\cap X>1$, and so $N\leq X$ by (3). Suppose that $H$ is $\Pi$-supplemented in $G$. Then $H$ has a supplement $T$ in $G$ such that $H\cap T\leq I\leq H$, where $I$ satisfies $\Pi$-property in $G$. Since $H\cap T=I\cap T$ and $N\leq T$, $L=H\cap N=I\cap N$. It follows that $|G:N_G(L)|=|G:N_G(I\cap N)|$ is a $p$-number. As $L\unlhd G_p$, we have that $L\unlhd G$. Then by (3), $L=1$, and so $|N|=p$, a contradiction. We may therefore, assume that $H$ has a $p$-supersoluble supplement $T$ in $G$. Let $F$ be the canonical local satellite of $\mathcal{U}_p$ such that $F(p)=\mathcal{G}_pF(p)=\mathcal{U}_p\cap \mathcal{G}_p\mathcal{A}(p-1)$, where $\mathcal{A}(p-1)$ denotes the class of finite abelian groups of exponent $p-1$ and $F(q)=\mathcal{U}_p$ for all primes $q\neq p$. By Lemma 2.4, $T/C_T(N)\in F(p)$. Since $P\leq C_G(N)$, we have that $G/C_G(N)\cong T/C_T(N)\in F(p)$. Then by Lemma 2.4 again, $N\leq Z_{\mathcal{U}_p}(G)$, and so $|N|=p$. The final contradiction ends the proof.\qed\par
\medskip
\noindent\textbf{Proof of Theorem C.} Suppose that the result is false and let $(G,E)$ be a counterexample for which $|G|+|E|$ is minimal. We proceed via the following steps.\par
(1) \textit{$O_{p'}(G)=1$.}\par
If not, by Lemma 2.1(2), $(G/O_{p'}(G), EO_{p'}(G)/O_{p'}(G))$ satisfies the hypothesis of this theorem. By the choice of $(G,E)$, we have that $EO_{p'}(G)/O_{p'}(G)$ is $p$-nilpotent, and so $E$ is $p$-nilpotent, a contradiction.\par
(2) \textit{$O_p(E)>1$.}\par
Suppose that $O_p(E)=1$ and let $N$ be a minimal normal subgroup of $G$ contained in $E$. Since $O_{p'}(G)=1$, $p\mid |N|$. Then we discuss three possible cases below:\par
(i) \textit{Case $1:$ $|N|_p<|D|$.}\par
In this case, by Lemma 2.2, $(G/N,E/N)$ satisfies the hypothesis of this theorem. By the choice of $(G,E)$, $E/N$ is $p$-nilpotent. Let $A/N$ be the normal $p$-complement of $E/N$. Then obviously, $A\unlhd G$ and $|A|_p=|N|_p<|D|$. By Lemma 2.2, $(G/A,E/A)$ satisfies the hypothesis of Theorem B. Therefore, $E/A\leq Z_\mathcal{U}(G/A)$. If $p|D|<|P|$, then we may take a normal subgroup $L$ of $G$ such that $A\leq L<E$ and $|L|_p=p|D|$. Clearly, $(G,L)$ satisfies the hypothesis of this theorem. Then by the choice of $(G,E)$, $L$ is $p$-nilpotent, and so $N$ is $p$-nilpotent. Since $O_{p'}(G)=1$, $N$ is a $p$-group. Hence $N\leq O_p(E)$, which is absurd.\par
Thus we have that $p|D|=|P|$. Then by hypothesis, every maximal subgroup of $P$ either is $\Pi$-supplemented in $G$ or has a $p$-supersoluble supplement in $G$. If every maximal subgroup of $P$ has a $p$-supersoluble supplement in $E$, then since $(|E|,p-1)=1$, every maximal subgroup of $P$ has a $p$-nilpotent supplement in $E$. By Lemmas 2.11 and 2.12, $E$ is $p$-nilpotent, a contradiction. Hence $P$ has a maximal subgroup $P_1$ such that $P_1$ is $\Pi$-supplemented in $G$ and $P_1$ does not have a $p$-supersoluble supplement in $E$. Then $P_1$ has a supplement $T$ in $G$ such that $T\cap E$ is not $p$-supersoluble and $P_1\cap T\leq I\leq P_1$, where $I$ satisfies $\Pi$-property in $G$. This implies that $|G:N_G(I\cap N)|$ is a $p$-number, and so $I\cap N\leq O_p(E)=1$. It follows that $P_1\cap T\cap N=I\cap N=1$. As $|T\cap E:P\cap T|=|E:P|$, $P\cap T$ is a Sylow $p$-subgroup of $T\cap E$. This induces that $P\cap T\cap N$ is a Sylow $p$-subgroup of $T\cap N$. Note that $|P\cap T\cap N|=|P\cap T\cap N:P_1\cap T\cap N|=|P_1(P\cap T\cap N):P_1|\leq p$. Hence $|T\cap N|_p\leq p$. Since $T\cap E/T\cap N\cong (TN\cap E)/N\leq E/N$ is $p$-nilpotent, by Lemma 2.13(2), $T\cap E$ is $p$-nilpotent, a contradiction.\par
(ii) \textit{Case $2:$ $|N|_p>|D|$.}\par
In this case, if $N<E$, then $(G,N)$ satisfies the hypothesis of this theorem. By the choice of $(G,E)$, $N$ is $p$-nilpotent. Since $O_{p'}(G)=1$, $N$ is a $p$-group, which is absurd. Hence $N=E$. By hypothesis, for every proper subgroup $H$ of $P$ with $|H|=p^n|D|$ ($n=0,1$) or 4 (when $|D|=1$ and $P$ is not quaternion-free), $H$ either is $\Pi$-supplemented in $G$ or has a $p$-supersoluble supplement in $G$. If $H$ is $\Pi$-supplemented in $G$, then $H$ has a supplement $T$ in $G$ such that $H\cap T\leq I\leq H$, where $I$ satisfies $\Pi$-property in $G$. It follows that $|G:N_G(I)|$ is a $p$-number, and so $I\leq O_p(E)=1$. Hence $H$ either is complemented in $G$ or has a $p$-supersoluble supplement in $G$. If $E<G$, then clearly, $H$ either is complemented in $E$ or has a $p$-supersoluble supplement in $E$. This shows that $(E,E)$ satisfies the hypothesis of this theorem. By the choice of $(G,E)$, $E$ is $p$-nilpotent, a contradiction. Thus $G=E$ is a non-abelian simple group. By Feit-Thompson Theorem, $p=2$.\par
If every maximal subgroup of $P$ has a $2$-supersoluble supplement in $G$, then by Lemmas 2.11 and 2.12, $G$ is $2$-nilpotent, which is impossible. This shows that $P$ has a maximal subgroup which does not have a $2$-supersoluble supplement in $G$. Suppose that $2|D|<|P|$. Then $P$ has subgroups $H_1$ and $H_2$ with $|H_1|=|D|$ and $|H_2|=2|D|$ such that $H_1$ and $H_2$ are complemented in $G$. Let $T_1$ and $T_2$ be complements of $H_1$ and $H_2$ in $G$, respectively. Then $|G:T_1|=2^{r}$ and $|G:T_2|=2^{r+1}$ such that $T_2$ is not a $2'$-Hall subgroup of $G$. If $T_1=G$, then $|G:T_2|=2$, and so $T_2\unlhd G$, which is absurd. Hence $r\geq 1$. If $r\leq 2$, then $G\cong G/(T_1)_G\apprle S_4$, where $S_4$ denotes the symmetric group of degree $4$, and so $G$ is soluble, a contradiction. Thus $r\geq 3$. By Lemma 2.14, $G$ is not a non-abelian simple group, which is impossible. Now assume that $2|D|=|P|$. Then $P$ has a maximal subgroup $H$ such that $H$ is complemented in $G$ such that every complement $T$ of $H$ in $G$ is not $2$-supersoluble. However, since $|T|_2=2$, $T$ is $2$-supersoluble, which contradicts our assumption.\par
(iii) \textit{Case $3:$ $|N|_p=|D|$.}\par
In this case, if $p|D|=|P|$, then $|E/N|_p=p$. Hence by Lemma 2.13(2), $E/N$ is $p$-nilpotent. With a similar argument as in the proof of Case 1 of (2), we can get a contradiction. Now assume that $p|D|<|P|$. Let $E/A$ be a chief factor of $G$ such that $N\leq A$. If $|A|_p>|N|_p=|D|$, then $(G,A)$ satisfies the hypothesis of this theorem. By the choice of $(G,E)$, $A$ is $p$-nilpotent. Since $O_{p'}(G)=1$, $A$ is a $p$-group, a contradiction. Hence $|A|_p=|N|_p=|D|$. By hypothesis, every subgroup of $P$ of order $|H|=p|D|$ either is $\Pi$-supplemented in $G$ or has a $p$-supersoluble supplement in $G$. Then by Lemma 2.2, every subgroup of $PA/A$ of order $p$ either is $\Pi$-supplemented in $G/A$ or has a $p$-supersoluble supplement in $G/A$.\par
Suppose that there exists a subgroup $H/A$ of $PA/A$ of order $p$ such that $H/A$ is $\Pi$-supplemented, but not complemented in $G/A$. Then clearly, $H/A$ satisfies $\Pi$-property in $G/A$. This implies that $|G/A:N_{G/A}(H/A)|$ is a $p$-number, and so $H/A\leq O_p(E/A)$. Hence $E/A=O_p(E/A)$ for $E/A$ is a chief factor of $G$. Consequently, $E/A$ is an elementary abelian $p$-group. Then $(G/A,E/A)$ satisfies the hypothesis of Theorem B. Thus $E/A\leq Z_\mathcal{U}(G/A)$. This induces that $|E/A|=p$, and so $p|D|=p|A|_p=|P|$, which is contrary to our assumption. Therefore, every subgroup of $PA/A$ of order $p$ either is complemented in $G/A$ or has a $p$-supersoluble supplement in $G/A$. Now we will show that $E/A$ is $p$-nilpotent. If $PA/A$ has a subgroup of order $p$ which has a $p$-supersoluble supplement in $G/A$, but is not complemented in $G/A$, then clearly, $G/A$ is $p$-supersoluble, and so is $E/A$. Since $(|E/A|,p-1)=1$, $E/A$ is $p$-nilpotent. Now consider that every subgroup of $PA/A$ of order $p$ is complemented in $G/A$. Then by Lemma 2.15, $E/A$ is also $p$-nilpotent. Since $p\mid|E/A|$, $E/A$ is an elementary abelian $p$-group. As discussed above, we can obtain that $|E/A|=p$, and thus $p|D|=p|A|_p=|P|$. The final contradiction shows that (2) holds.\par
(3) \textit{Final contradiction.}\par
Since $O_p(E)>1$, let $N$ be a minimal normal subgroup of $G$ contained in $O_p(E)$. Then we discuss three possible cases as follows:\par
(i) \textit{Case $1:$ $|N|<|D|$.}\par
In this case, by Lemma 2.1(2), $(G/N,E/N)$ satisfies the hypothesis of this theorem. By the choice of $(G,E)$, $E/N$ is $p$-nilpotent. Let $A/N$ be the normal $p$-complement of $E/N$. Since $|A|_p=|N|<|D|$,  by Lemma 2.2, $(G/A, E/A)$ satisfies the hypothesis of Theorem B, and so $E/A\leq Z_\mathcal{U}(G/A)$. If $p|D|<|P|$, then we may take a normal subgroup $L$ of $G$ such that $A\leq L<E$ and $|L|_p=p|D|$. It is easy to see that $(G,L)$ satisfies the hypothesis of this theorem. Then by the choice of $(G,E)$, $L$ is $p$-nilpotent. Since $O_{p'}(G)=1$, $L$ is a $p$-group. It follows that $E$ is a $p$-group, a contradiction.\par
We may, therefore, assume that $p|D|=|P|$. By hypothesis, every maximal subgroup of $P$ either is $\Pi$-supplemented in $G$ or has a $p$-supersoluble supplement in $G$. Since $E/N$ is $p$-nilpotent, by Lemma 2.16, $N\nleq \Phi(G)$. Thus $G$ has a maximal subgroup $M$ such that $N\nleq M$. Obviously, $M\cap N=1$. As $E/N$ is $p$-nilpotent, $M\cap E$ is $p$-nilpotent. Let $G_p$ be a Sylow $p$-subgroup of $G$ and $G_{p_1}$ a maximal subgroup of $G_p$ containing $G_p\cap M$. Then $G_p=G_{p_1}N$. Let $P_1=G_{p_1}\cap P$. Since $|P:P_1|=|G_p:G_{p_1}|=p$, $P_1$ is a maximal subgroup of $P$ such that $P=P_1N$. If $P_1$ is $\Pi$-supplemented in $G$, then $P_1$ has a supplement $T$ in $G$ such that $P_1\cap T\leq I\leq P_1$, where $I$ satisfies $\Pi$-property in $G$. It follows that $|G:N_G(I\cap N)|$ is a $p$-number. If $I\cap N>1$, then $N=(I\cap N)^G=(I\cap N)^{G_p}\leq G_{p_1}$, a contradiction. Hence $I\cap N=1$, and thus $P_1\cap T\cap N=1$. Since $T\cap E/T\cap N\cong (TN\cap E)/N$ is $p$-nilpotent and $|T\cap N|=|T\cap N:P_1\cap T\cap N|=|P_1(T\cap N):P_1|\leq p$, $T\cap E$ is $p$-nilpotent by Lemma 2.13(2). Consequently, no matter $P_1$ is $\Pi$-supplemented in $G$ or has a $p$-supersoluble supplement in $G$, $P_1$ has a $p$-nilpotent supplement $T_1$ in $E$ for $(|E|,p-1)=1$. Let $(M\cap E)_{p'}$ and ${(T_1)}_{p'}$ be the normal $p$-complements of $M\cap E$ and $T_1$, respectively. Then $(M\cap E)_{p'}$ and ${(T_1)}_{p'}$ are $p'$-Hall subgroups of $E$. By Lemma 2.11, $E$ is a $C_{p'}$-group. This implies that $E$ has an element $g$ such that ${(T_1)}_{p'}^g=(M\cap E)_{p'}$. Considering the fact that $T_1\leq N_E({(T_1)}_{p'})$, we may let $g\in P_1$. It follows that $E=P_1N_E({(T_1)}_{p'})^g=P_1N_E((M\cap E)_{p'})$. Since $O_{p'}(G)=1$ and $M\leq N_G((M\cap E)_{p'})$, we have that $N_G((M\cap E)_{p'})=M$. This implies that $E=P_1(M\cap E)$. As $P\cap M\leq G_{p_1}\cap P=P_1$, $P=P_1(P\cap M)=P_1$, which is impossible.\par
(ii) \textit{Case $2:$ $|N|>|D|$.}\par
In this case, $(G,N)$ satisfies the hypothesis of Theorem B. Hence $N\leq Z_\mathcal{U}(G)$, and so $|N|=p$. It follows that $|D|=1$. As $(G,O_p(E))$ satisfies the hypothesis of Theorem B, $O_p(E)\leq Z_\mathcal{U}(G)\leq Z_\mathcal{U}(E)$. Since $(|E|,p-1)=1$, it is easy to see that $O_p(E)\leq Z_\infty(E)$. Let $A/O_p(E)$ be a chief factor of $G$ below $E$. If $A<E$, then $(G,A)$ satisfies the hypothesis of this theorem. By the choice of $(G,E)$, $A$ is $p$-nilpotent. Then since $O_{p'}(G)=1$, $A$ is a $p$-group. This shows that $A\leq O_p(E)$, which is absurd. Hence $E/O_p(E)$ is a chief factor of $G$. If $p\nmid|E/O_p(E)|$, then $E$ is $p$-nilpotent for $O_p(E)\leq Z_\infty(E)$, a contradiction. Thus $p\mid|E/O_p(E)|$. Obviously, $E$ is not soluble. Thus by Feit-Thompson Theorem, $p=2$.\par
Now let $V$ be a minimal non-$2$-nilpotent group contained in $E$. By \cite[Chapter IV, Satz 5.4]{Hup1}, $V$ is a minimal non-nilpotent group such that $V=V_2\rtimes V_q$, where $V_2$ is the Sylow 2-subgroup of $V$ and $V_q$ is a Sylow $q$-subgroup of $V$ with $q>2$. Without loss of generality, we may let $V_2\leq P$. Then by \cite[Chapter VII, Theorem 6.18]{Doe}, $V_2/\Phi(V_2)$ is a $V$-chief factor; $\Phi(V)=Z_\infty(V)$; $\Phi(V_2)=V_2\cap \Phi(V)$; and $V_2$ has exponent $2$ or $4$ (when $V_2$ is non-abelian). It follows that $O_2(E)\cap V_2\leq Z_\infty(E)\cap V_2\leq Z_\infty(V)\cap V_2=\Phi(V)\cap V_2=\Phi(V_2)$. Therefore, $V_2$ has an element $x$ which is not contained in $O_2(E)$. Let $H=\langle x \rangle$. Then $|H|=2$ or $4$ (when $V_2$ is non-abelian). If $V_2$ is non-abelian and quaternion-free, then $V_2$ has a characteristic subgroup of index 2 by Lemma 2.9. This implies that $|V_2/\Phi(V_2)|=2$, and so $V_2$ is cyclic, which contradicts our assumption. Therefore, $|H|=2$ or $4$ (when $V_2$ is not quaternion-free). By hypothesis, $H$ either is $\Pi$-supplemented in $G$ or has a $p$-supersoluble supplement in $G$. 
Let $X$ be any supplement of $H$ in $G$. Suppose that $X<G$. Then $G/X_G\apprle S_4$ for $|G:X|\leq 4$, where $S_4$ denotes the symmetric group of degree $4$. Thus $E/X_G\cap E$ is soluble. 
Since $X_G\cap E<E$ and $(G,X_G\cap E)$ satisfies the hypothesis of this theorem, $X_G\cap E$ is $2$-nilpotent by the choice of $(G,E)$. This induces that $E$ is soluble, which is impossible. Therefore, $G$ is the unique supplement of $H$ in $G$. Since $G$ is not $2$-supersoluble, $H$ is $\Pi$-supplemented in $G$, and so $H$ satisfies $\Pi$-property in $G$. Then $|G:N_G(HO_2(E))|$ is a $2$-number. This implies that $H\leq O_2(E)$, a final contradiction of (ii).\par
(iii) \textit{Case $3:$ $|N|=|D|$.}\par
In this case, if $p|D|=|P|$, then $|E/N|_p=p$, and thus $E/N$ is $p$-nilpotent by Lemma 2.13(2). With a similar discussion as in the proof of Case 1 of (3), we can get a contradiction. Hence $p|D|<|P|$. Let $E/A$ be a chief factor of $G$ such that $N\leq A$. If $|A|_p=|N|=|D|$, then a contradiction can be derived in a similar way as in Case 3 of (2). Now we may assume that $|A|_p>|N|=|D|$. Then $(G,A)$ satisfies the hypothesis of this theorem. By the choice of $(G,E)$, $A$ is $p$-nilpotent. Since $O_{p'}(G)=1$, $A$ is a $p$-group. It follows that $(G,A)$ satisfies the hypothesis of Theorem B. Hence $A\leq Z_\mathcal{U}(G)$, and so $|N|=|D|=p$. This case can be viewed as a special case of Case 2 of (3) (we may take $|N|=p$ and $|D|=1$), and this fact yields a contradiction. The theorem is thus proved. \qed \par
\medskip
\noindent\textbf{Proof of Theorem A.} Let $p$ be the smallest prime divisor of $|X|$ and $P$ a Sylow $p$-subgroup of $X$. If $P$ is cyclic, then by Lemma 2.13(1), $X$ is $p$-nilpotent. Now assume that $P$ is not cyclic. Then by Theorem C, $X$ is also $p$-nilpotent. Let $X_{p'}$ be the normal $p$-complement of $X$. Then $X_{p'}\unlhd G$. If $P$ is cyclic, then $X/X_{p'}\leq Z_\mathcal{U}(G/X_{p'})$. Now consider that $P$ is not cyclic. Then by Lemma 2.1(2), $(G/X_{p'},X/X_{p'})$ satisfies the hypothesis of Theorem B. Hence we also have that $X/X_{p'}\leq Z_\mathcal{U}(G/X_{p'})$.\par
Let $q$ be the second smallest prime divisor of $|X|$ and $Q$ a Sylow $q$-subgroup of $X$. With a similar argument as above, we can get that $X_{p'}$ is $q$-nilpotent and $X_{p'}/X_{\{p,q\}'}\leq Z_\mathcal{U}(G/X_{\{p,q\}'})$, where $X_{\{p,q\}'}$ is the normal $q$-complement of $X_{p'}$. The rest may be deduced by analogy. Therefore, we obtain that $X\leq Z_\mathcal{U}(G)\leq Z_\mathcal{F}(G)$. It follows from Lemma 2.17 that $E\leq Z_\mathcal{F}(G)$. Then by Lemma 2.18, $G\in \mathcal{F}$ as desired.\qed\par
\section{Further Applications}
In this section, we will show that the subgroups of $G$ which satisfy a certain known embedding property mentioned below are all $\Pi$-supplemented in $G$. For the sake of simplicity, we only focus on most recent embedding properties.\par
Recall that a subgroup $H$ of $G$ is called to be a \textit{CAP-subgroup} if $H$ either covers or avoids every chief factor of $G$. Let $\mathcal{F}$ be a saturated formation. A subgroup $H$ of $G$ is said to be \textit{$\mathcal{F}$-hypercentrally embedded} \cite{Ezq} in $G$ if $H^G/H_G\leq Z_{\mathcal{F}}(G/H_G)$. A subgroup $H$ of $G$ is called to be \textit{S-quasinormal} (or S-permutable) in $G$ if $H$ permutes with every Sylow subgroup of $G$. A subgroup $H$ of $G$ is said to be \textit{S-semipermutable} \cite{Che} in $G$ if $H$ permutes with every Sylow $p$-subgroup of $G$ such that $(p,|H|)=1$. A subgroup $H$ of $G$ is called to be \textit{S-quasinormally embedded} \cite{Bal7} in $G$ if every Sylow subgroup of $H$ is a Sylow subgroup of some S-quasinormal subgroup of $G$. A subgroup $H$ of $G$ is said to be \textit{S-conditionally permutable} \cite{Hua} in $G$ if $H$ permutes with at least one Sylow $p$-subgroup of $G$ for every $p\in\pi(G)$.\par
\medskip
\noindent\textbf{Proposition 4.1.} Let $H$ be a subgroup of $G$. Then $H$ satisfies $\Pi$-property in $G$ if one of the following holds:\par
(1) $H$ is a CAP-subgroup of $G$.\par
(2) $H$ is $\mathcal{U}$-hypercentrally embedded in $G$.\par
(3) $H$ is S-quasinormal in $G$.\par
(4) $H$ is a $p$-group and $H$ is S-semipermutable in $G$.\par
(5) $H^G$ is soluble and $H$ is S-quasinormally embedded in $G$.\par
(6) $H^G$ is soluble and $H$ is S-conditionally permutable in $G$.\par
\smallskip
\noindent\textit{Proof.} Statements (1)-(3) and (5)-(6) were proved in \cite{Li1}, and the proof of \cite[Proposition 2.4]{Li1} still works for statement (4). \par
\medskip
Recall that a subgroup $H$ of $G$ is called to be a \textit{CAS-subgroup} \cite{Wei1} if there exists a subgroup $T$ of $G$ such that $G=HT$ and $H\cap T$ is a CAP-subgroup of $G$. Let $\mathcal{F}$ be a saturated formation. A subgroup $H$ of $G$ is said to be \textit{$\mathcal{F}$-supplemented} \cite{Guo7} in $G$ if there exists a subgroup $T$ of $G$ such that $G=HT$ and $(H\cap T)H_G/H_G\leq Z_{\mathcal{F}}(G/H_G)$. A subgroup $H$ of $G$ is called to be \textit{weakly $s$-supplemented} \cite{Ski3} in $G$ if there exists a subgroup $T$ of $G$ such that $G=HT$ and $H\cap T\leq H_{sG}$, where $H_{sG}$ denotes the subgroup generated by all those subgroups of $H$ which are S-quasinormal in $G$. A subgroup $H$ of $G$ is said to be \textit{weakly $\bar{s}$-supplemented} \cite{Xu} in $G$ if there exists a subgroup $T$ of $G$ such that $G=HT$ and $H\cap T\leq H_{\bar{s}G}$, where $H_{\bar{s}G}$ denotes the subgroup generated by all those subgroups of $H$ which are S-semipermutable in $G$. A subgroup $H$ of $G$ is called to be \textit{weakly $s$-supplementally embedded} \cite{Zha} in $G$ if there exists a subgroup $T$ of $G$ such that $G=HT$ and $H\cap T\leq H_{se}$, where $H_{se}$ denotes an S-quasinormally embedded subgroup of $G$ contained in $H$. A subgroup $H$ of $G$ is said to be \textit{completely $c$-permutable} \cite{Guo9} in $G$ if for every subgroup $T$ of $G$, there exists some $x\in \langle H,T\rangle$ such that $HT^x=T^xH$. A subgroup $H$ of $G$ is called to be \textit{weakly $c$-permutable} \cite{Guo8} in $G$ if there exists a subgroup $T$ of $G$ such that $G=HT$ and $H\cap T$ is completely $c$-permutable in $G$.\par
\medskip
\noindent\textbf{Proposition 4.2.} Let $H$ be a subgroup of $G$. Then $H$ is $\Pi$-supplemented in $G$ if one of the following holds:\par
(1) $H$ is a CAS-subgroup of $G$.\par
(2) $H$ is $\mathcal{U}$-supplemented in $G$.\par
(3) $H$ is weakly $s$-supplemented in $G$.\par
(4) $H$ is a $p$-group and $H$ is weakly $\bar{s}$-supplemented in $G$.\par
(5) $H^G$ is soluble and $H$ is weakly $s$-supplementally embedded in $G$.\par
(6) $H^G$ is soluble and $H$ is weakly $c$-permutable in $G$.\par
\smallskip
\noindent\textit{Proof.} Note that by \cite[Satz 2]{Keg}, $H_{sG}$ is S-quasinormal in $G$, and if $H$ is a $p$-group, then by definition, $H_{\bar{s}G}$ is S-semipermutable in $G$. Also, every completely $c$-permutable subgroup of $G$ is clearly S-conditionally permutable in $G$. Then Proposition 4.2 directly follows from Proposition 4.1.\par
\medskip
By the above proposition, a large number of previous results are immediate consequences of our theorems. We omit further details, and readers may refer to the relevant literature for more information.\par
\bibliographystyle{plain}
\bibliography{expbib}

\begin{thebibliography}{10}

\bibitem{Asa1}
M.~Asaad.
\newblock Finite groups with certain subgroups of {S}ylow subgroups
  complemented.
\newblock {\em J. Algebra}, 323:1958--1965, 2010.

\bibitem{Asa2}
M.~Asaad.
\newblock On $c$-supplemented subgroups of finite groups.
\newblock {\em J. Algebra}, 362:1--11, 2012.

\bibitem{Bal7}
A.~Ballester-Bolinches and M.~C. Pedraza-Aguilera.
\newblock Sufficient conditions for supersolubility of finite groups.
\newblock {\em J. Pure Appl. Algebra}, 127:113--118, 1998.

\bibitem{Bal3}
A.~Ballester-Bolinches and M.~D. P\'erez-Ramos.
\newblock On $\mathfrak{F}$-subnormal subgroups and {F}rattini-like subgroups
  of a finite group.
\newblock {\em Glasgow Math. J.}, 36:241--247, 1994.

\bibitem{Bal2}
A.~Ballester-Bolinches, Y.~Wang, and X.~Guo.
\newblock $c$-supplemented subgroups of finite groups.
\newblock {\em Glasgow Math. J.}, 42:383--389, 2000.

\bibitem{Che2}
X.~Chen, W.~Guo, and A.~N. Skiba.
\newblock Some conditions under which a finite group belongs to a {B}aer-local
  formation.
\newblock {\em Comm. Algebra}, 2013.
\newblock (accepted).

\bibitem{Che}
Z.~Chen.
\newblock On a theorem of {S}rinivasan (in {C}hinese).
\newblock {\em J. Southwest Normal Univ. Nat. Sci.}, 12(1):1--4, 1987.

\bibitem{Cro}
F.~Cross.
\newblock Conjugacy of odd order {H}all subgroups.
\newblock {\em Bull. London Math. Soc.}, 19:311--319, 1987.

\bibitem{Doe}
K.~Doerk and T.~Hawkes.
\newblock {\em Finite Soluble Groups}.
\newblock Walter de Gruyter, Berlin/New York, 1992.

\bibitem{Dor}
L.~Dornhoff.
\newblock M-groups and 2-groups.
\newblock {\em Math. Z.}, 100:226--256, 1967.

\bibitem{Ezq}
L.~M. Ezquerro and X.~Soler-Escriv\`{a}.
\newblock Some permutability properties related to $\mathcal{F}$-hypercentrally
  embedded subgroups of finite groups.
\newblock {\em J. Algebra}, 264:279--295, 2003.

\bibitem{Gag}
T.~M. Gagen.
\newblock {\em Topics in Finite Groups}.
\newblock Cambridge, London/New York/Melbourne, 1976.

\bibitem{Gor}
D.~Gorenstein.
\newblock {\em Finite Groups}.
\newblock Harper \& Row Publishers, New York/Evanston/London, 1968.

\bibitem{Guo1}
W.~Guo.
\newblock {\em The Theory of Classes of Groups}.
\newblock Kluwer, Dordrecht, 2000.

\bibitem{Guo7}
W.~Guo.
\newblock On $\mathcal{F}$-supplemented subgroups of finite groups.
\newblock {\em Manuscripta Math.}, 127:139--150, 2008.

\bibitem{Guo8}
W.~Guo and S.~Chen.
\newblock Weakly $c$-permutable subgroups of finite groups.
\newblock {\em J. Algebra}, 324:2369--2381, 2010.

\bibitem{Guo9}
W.~Guo, K.~P. Shum, and A.~N. Skiba.
\newblock Conditionally permutable subgroups and supersolubility of finite
  groups.
\newblock {\em Southeast Asian Bull. Math.}, 29:493--510, 2005.

\bibitem{Guo}
W.~Guo and A.~N. Skiba.
\newblock On {$\mathcal{F}\mathnormal{\Phi}^*$}-hypercentral subgroups of
  finite groups.
\newblock {\em J. Algebra}, 372:275--292, 2012.

\bibitem{Hua}
J.~Huang and W.~Guo.
\newblock S-conditionally permutable subgroups of finite groups (in {C}hinese).
\newblock {\em Chin. Ann. Math. Ser. A}, 28(1):17--26, 2007.

\bibitem{Hup1}
B.~Huppert.
\newblock {\em Endliche Gruppen I}.
\newblock Springer-Verlag, 1968.

\bibitem{Hup}
B.~Huppert and N.~Blackburn.
\newblock {\em Finite groups III}.
\newblock Springer-Verlag, Berlin/Heidelberg, 1982.

\bibitem{Keg}
O.~H. Kegel.
\newblock Sylow-gruppen und subnormalteiler endlicher gruppen.
\newblock {\em Math. Z.}, 78:205--221, 1962.

\bibitem{Li1}
B.~Li.
\newblock On {$\Pi$}-property and {$\Pi$}-normality of subgroups of finite
  groups.
\newblock {\em J. Algebra}, 334:321--337, 2011.

\bibitem{Li6}
B.~Li.
\newblock Finite groups with {$\Pi$}-supplemented minimal subgroups.
\newblock {\em Comm. Algebra}, 41:2060--2070, 2013.

\bibitem{Ski3}
A.~N. Skiba.
\newblock On weakly $s$-permutable subgroups of finite groups.
\newblock {\em J. Algebra}, 315:192--209, 2007.

\bibitem{Ski}
A.~N. Skiba.
\newblock On two questions of {L. A. S}hemetkov concerning hypercyclically
  embedded subgroups of finite groups.
\newblock {\em J. Group Theory}, 13:841--850, 2010.

\bibitem{Ski1}
A.~N. Skiba and L.~A. Shemetkov.
\newblock Multiply $\mathcal{L}$-composition formations of finite groups.
\newblock {\em Ukr. Math. J.}, 52(6):898--913, 2000.

\bibitem{War}
H.~N. Ward.
\newblock Automorphisms of quaternion-free 2-groups.
\newblock {\em Math. Z.}, 112:52--58, 1969.

\bibitem{Wei1}
H.~Wei and Y.~Wang.
\newblock On {CAS}-subgroups of finite groups.
\newblock {\em Israel J. Math.}, 159:175--188, 2007.

\bibitem{Xu}
Y.~Xu and X.~Li.
\newblock Weakly $s$-semipermutable subgroups of finite groups.
\newblock {\em Front. Math. China}, 6(1):161--175, 2011.

\bibitem{Zha}
T.~Zhao, X.~Li, and Y.~Xu.
\newblock On weakly $s$-supplemently embedded subgroups of finite groups.
\newblock {\em Proc. {E}dinburgh Math. Soc.}, 54:799--807, 2011.

\end{thebibliography}
\end{document}